\documentclass[11pt]{amsart}

\usepackage{fullpage}

\usepackage{amssymb,pb-diagram,lamsarrow,pb-lams}

\usepackage{rotating}

\usepackage{array}
\newcolumntype{C}{>{$}c<{$}}
\newcolumntype{R}{>{$}r<{$}}
\newcolumntype{L}{>{$}l<{$}}
\newcolumntype{S}{>{$\scriptstyle}l<{$}}

\def\to{\mathchoice
{\longrightarrow}
{\rightarrow}
{\rightarrow}
{\rightarrow}}

\def\mapsto{\DOTSB\mapstochar\to}

\def\hookrightarrow{\mathchoice
{\DOTSB\lhook\joinrel\relbar\joinrel\rightarrow}
{\DOTSB\lhook\joinrel\rightarrow}
{\DOTSB\lhook\joinrel\rightarrow}
{\DOTSB\lhook\joinrel\rightarrow}}

\newtheorem{theorem}{Theorem}
\newtheorem*{theorem*}{Theorem}
\newtheorem{proposition}[theorem]{Proposition}
\newtheorem{lemma}[theorem]{Lemma}
\newtheorem{corollary}[theorem]{Corollary}

\theoremstyle{definition}

\newcommand{\Z}{\mathbb{Z}}
\newcommand{\C}{\mathbb{C}}

\newcommand{\Q}{\mathbb{Q}}

\newcommand{\om}{\omega}
\newcommand{\eps}{\varepsilon}
\renewcommand{\o}{\otimes}
\newcommand{\p}{\partial}
\renewcommand{\*}{\cdot}

\newcommand{\<}{\langle}
\renewcommand{\>}{\rangle}

\renewcommand{\[}{{[\![}}

\newcommand{\CM}{\mathcal{M}}
\newcommand{\CO}{\mathcal{O}}
\newcommand{\I}{\mathcal{I}}

\newcommand{\CH}{\mathcal{H}}
\newcommand{\Mbar}{\overline{\mathcal{M}}}

\renewcommand{\SS}{\mathbb{S}}
\DeclareMathOperator{\SP}{Sp}

\newcommand{\Om}{\Omega}

\DeclareMathOperator{\SL}{SL}

\DeclareMathOperator{\Hom}{Hom}

\newcommand{\CC}{\mathsf{C}}
\newcommand{\QQ}{\mathsf{Q}}
\newcommand{\TT}{\mathsf{T}}
\newcommand{\OO}{\mathsf{O}}
\newcommand{\II}{\mathsf{I}}

\newcommand{\V}{\mathsf{V}}
\newcommand{\W}{\mathsf{W}}
\newcommand{\VV}{\mathbb{V}}
\newcommand{\WW}{\mathbb{W}}
\newcommand{\Schur}{S}

\newcommand{\PP}{\mathbb{P}}

\newcommand{\Dual}{\vee}
\newcommand{\Id}{I}
\DeclareMathOperator{\Aut}{Aut}

\DeclareMathOperator{\Wedge}{\Lambda}
\DeclareMathOperator{\ch}{ch}

\DeclareMathOperator{\Spec}{Spec}

\begin{document}

\title{Euler characteristics of local systems on $\CM_2$}

\author{E. Getzler}

\address{Department of Mathematics, Northwestern University, Evanston, IL
60208-2730}

\email{getzler@math.nwu.edu}


\thanks{This paper was completed while the author was a guest of the School
of Mathematical Sciences of the Australian National University}
\thanks{The author is partially supported under NSF grant DMS-9704320}

\begin{abstract}
We calculate the Euler characteristics of the local systems
$\Schur^k\VV\o\Schur^\ell\Wedge^2\VV$ on the moduli space $\CM_2$ of curves
of genus $2$, where $\VV$ is the rank $4$ local system $R^1\pi_*\C$.
\end{abstract}

\maketitle

\section{Introduction}

Let $\pi_g:\CM_{g,1}\to\CM_g$ be the universal curve of genus $g$,
$g\ge2$. The local system $\VV = R^1\pi_{g*}\C$ is symplectic of
rank $2g$. Given a sequence $(k_1,\dots,k_g)$ of nonnegative integers,
there is an associated local system
$$
\WW(1^{k_1}\dots g^{k_g}) = \Schur^{k_1}\bigl(\VV\bigr) \o
\Schur^{k_2}\bigl(\Wedge^2\VV\bigr) \o \dots \o
\Schur^{k_g}\bigl(\Wedge^g\VV\bigr)
$$
on $\CM_g$, with Euler characteristic $e_g(1^{k_1}\dots g^{k_g})=e\bigl(
\CM_g,\WW(1^{k_1}\dots g^{k_g}) \bigr)$, and generating function
$$
f_g(u_1,\dots,u_g) = \sum_{k_1,\dots,k_g=0}^\infty u_1^{k_1} \dots
u_g^{k_g} e_g(1^{k_1}\dots g^{k_g}) .
$$
In this paper, we calculate $f_2(u_1,u_2)$.

Let $\Gamma_g=\pi_1(\CM_g)$ be the genus $g$ mapping class group. There is
a homomorphism $\rho_g$ from $\Gamma_g$ to the symplectic group
$\SP(2g,\C)$, obtained by composition of the quotient map to $\SP(2g,\Z)$
with the inclusion $\SP(2g,\Z)\hookrightarrow\SP(2g,\C)$, and the Euler
characteristics $e_g(1^{k_1}\dots g^{k_g})$ may also be realized as the
Euler characteristic in group cohomology
$$
\sum_i (-1)^i \dim H^i\bigl( \Gamma_g , \rho^* \bigl( \Schur^{k_1}(\C^{2g})
\o \Schur^{k_2}(\Wedge^2\C^{2g}) \o \dots \o \Schur^{k_g}(\Wedge^g\C^{2g})
\bigr) \bigr) .
$$

We may illustrate our method by considering the analogous problem in genus
$1$; here, we must replace the universal curve by the fibration
$\pi_1:\CM_{1,2}\to\CM_{1,1}$. The generating function in this case equals
$$
f_1(u) = \sum_{k=0}^\infty u^k e\bigl( \CM_{1,1},\Schur^k\VV \bigr) =
\sum_{k=0}^\infty u^k e\bigl( \SL(2,\Z) ,\Schur^k\C^2 \bigr) .
$$
To calculate $f_1$, we stratify the coarse moduli space $|\CM_{1,1}|$ of
$\CM_{1,1}$, in other words the $j$-line, according to the automorphism
group of the elliptic curve $E(j)$; let $\CM_{1,1}(\Gamma)$ be the
subvariety of $|\CM_{1,1}|$ where $E(j)$ has automorphism group isomorphic
to $\Gamma$. There are three strata:
\begin{align*}
\CM_{1,1}(\CC_2) &= \C \setminus \{0,1728\} , \\
\CM_{1,1}(\CC_4) &= \{j=1728\} , \\
\CM_{1,1}(\CC_6) &= \{j=0\} .
\end{align*}
Denote the projection from the stack $\CM_{1,1}$ to $|\CM_{1,1}|$ by $\mu$.
If $\WW$ is a local system on $\CM_{1,1}$, we have
$$
e(\CM_{1,1},\WW) = e(|\CM_{1,1}|,\mu_*\WW) ,
$$
since $R^i\mu_*\WW=0$ for $i>0$. The Euler characteristic of a local
systems on a stratified space is the sum of the Euler characteristics over
the strata:
$$
e(|\CM_{1,1}|,\mu_*\WW) = e(\CM_{1,1}(\CC_2),\mu_*\WW) +
e(\CM_{1,1}(\CC_4),\mu_*\WW) + e(\CM_{1,1}(\CC_6),\mu_*\WW) .
$$
The restriction of the constructible sheaf $\mu_*\WW$ to a stratum
$\CM_{1,1}(\Gamma)$ is a local system, and hence its Euler characteristic
on this stratum is equal to the product of the Euler characteristic of
$\CM_{1,1}(\Gamma)$ and the rank of $\mu_*\WW$ restricted to
$\CM_{1,1}(\Gamma)$. (It is the failure of the analogous property for
stacks which necessitates the descent to the coarse moduli space
$|\CM_{1,1}|$.) We conclude that
\begin{multline*}
f_1(u) = \sum_{k=0}^\infty u^k \Bigl\{
e\bigl(\CM_{1,1}(\CC_2),\Schur^k\VV\bigr) +
e\bigl(\CM_{1,1}(\CC_4),\Schur^k\VV\bigr) +
e\bigl(\CM_{1,1}(\CC_6),\Schur^k\VV\bigr) \Bigr\} \\
= \sum_{k=0}^\infty u^k \Bigl\{
e\bigl(\CM_{1,1}(\CC_2)\bigr) \dim(\Schur^k\C^2)^{\CC_2} +
e\bigl(\CM_{1,1}(\CC_4)\bigr) \dim(\Schur^k\C^2)^{\CC_4} +
e\bigl(\CM_{1,1}(\CC_6)\bigr) \dim(\Schur^k\C^2)^{\CC_6} \Bigr\} .
\end{multline*}

The cyclic group $\CC_n$ is conjugate to the subgroup $\Bigl\{
\bigl(\begin{smallmatrix}z&0\\0&z^{-1}\end{smallmatrix}\bigr) \Bigm| z^n=1
\Bigr\}$ of $\SL(2,\C)$; it follows that
$$
\sum_{k=0}^\infty u^k e(\CC_n,\Schur^k\C^2) = \sum_{k=0}^\infty u^k
\dim(\Schur^k\C^2)^{\CC_n} = \frac{1+u^n}{(1-u^2)(1-u^n)} .
$$
It follows that
$$
f_1(u) = - \frac{1+u^2}{(1-u^2)^2} + \frac{1+u^4}{(1-u^2)(1-u^4)} +
\frac{1+u^6}{(1-u^2)(1-u^6)} = \frac{1-u^2-2u^4-u^6+u^8}{(1-u^4)(1-u^6)} .
$$

Our calculation in genus $2$ proceeds analogously: we use Bolza's
stratification of the coarse moduli space $|\CM_2|$ by the automorphism
group of the corresponding curve \cite{Bolza}. Denote the stratum
associated to the automorphism group $\Gamma$ by $\CM_2(\Gamma)$. The
contribution of each stratum to $e(\CM_2,\WW(1^k2^\ell))$ must be
calculated using the character theory of $\Gamma$; since $\Gamma$ is a
finite subgroup of $\SL(2,\C)$, the McKay correspondence allows this to be
done in terms of the associated Dynkin diagram. The only tricky point is
the calculation of the Euler characteristics of $\CM_2(\Gamma)$. The
hardest case is the affine surface $\CM_2(\CC_4)$; we prove in Section~3
that $e(\CM_2(\CC_4))=3$.

The original motivation for this work was the desire to calculate the
$\SS_n$-equivariant Euler characteristics of the moduli spaces
$\CM_{2,n}$. We explain how this may be done in Section 5.

Throughout this paper, $\eps_n$ denotes a primitive $n$th root of unity.
All varieties we consider are defined over the field of complex numbers
$\C$.

\section{Finite subgroups of $\SL(2,\C)$ and the McKay correspondence}

Given positive integers $p\ge q\ge r\ge2$, let $\<p,q,r\>$ be the group
with presentation
$$
\< S,T,U \mid S^p=T^q=U^r=STU \> .
$$
The element $STU$ is a central involution, which we denote by $-\Id$.

If $p^{-1}+q^{-1}+r^{-1}>1$, the group $\<p,q,r\>$ is finite, and its order
equals $4/(p^{-1}+q^{-1}+r^{-1}-1)$. This happens in the following cases.
$$\setlength{\extrarowheight}{4pt}
\begin{tabular}{|C|C|C|c|C|C|} \hline
(p,q,r) & \<p,q,r\> & \text{order} & Name & S \\[3pt] \hline
(n,2,2) & \QQ_{4n} & 4n & quaternionic, $n\ge2$ & \left( \begin{smallmatrix}
\eps_{2n} & 0 \\ 0 & \eps_{2n}^{-1} \end{smallmatrix} \right) \\[7pt]
(3,3,2) & \TT & 24 & binary tetrahedral & \frac{1}{\sqrt{2}} \left(
\begin{smallmatrix} \eps_8^{-1} & \eps_8^3 \\ \eps_8 & \eps_8
\end{smallmatrix} \right) \\[7pt]
(4,3,2) & \OO & 48 & binary octahedral & - \frac{1}{\sqrt{2}} \left(
\begin{smallmatrix} 1 & \eps_8 \\ \eps_8^3 & 1 \end{smallmatrix} \right)
\\[7pt]
(5,3,2) & \II & 120 & binary icosahedral & \frac{1}{\sqrt{5}} \left(
\begin{smallmatrix} \eps^3_5-1 & \eps_5^3-\eps_5 \\ \eps_5^4-\eps_5^2 &
\eps_5-1 \end{smallmatrix} \right) \\[8pt] \hline
\end{tabular}
\setlength{\extrarowheight}{0pt}$$

According to Klein \cite{Klein}, the non-abelian finite subgroups $\Gamma$
of $\SL(2,\C)$ are isomorphic to the finite groups $\<p,q,r\>$; any such
subgroup of $\SL(2,\C)$ is conjugate to the subgroup generated by the
element $S$ of $\SL(2,\C)$ listed in the table, and the element $U=\left(
\begin{smallmatrix} 0 & 1 \\ -1 & 0 \end{smallmatrix} \right)$.

The abelian subgroups of $\SL(2,\C)$ are all cyclic, and any such subgroup
of $\SL(2,\C)$ is conjugate to the subgroup generated by the element $T =
\Bigl( \begin{smallmatrix} \eps_n & 0 \\ 0 & \eps_n^{-1}
\end{smallmatrix} \Bigr)$ of $\SL(2,\C)$

We refer to the finite subgroups of $\SL(2,\C)$ as the Kleinian groups.  If
$\Gamma$ is such a group, let $\V$ be the two-dimensional fundamental
representation of $\Gamma$ induced by the embedding of $\Gamma$ in
$\SL(2,\C)$, and let $\V_k$ be the $k$th symmetric power $\Schur^k\V$ of
$\V$ (isomorphic to the space of binary forms of degree $k$). If $\Gamma$
is a Kleinian group $\Gamma$ containing $-\Id$ and $\W$ is an irreducible
representation of $\Gamma$, we call $\W$ even if $-\Id$ acts by $+1$ and
odd if its acts by $-1$; by Schur's lemma, these are the only
possibilities. For example, the fundamental representation $\V$ is odd.

There is a beautiful relationship between the character theory of Kleinian
groups and Dynkin diagrams, known as the McKay correspondence. If $\Gamma$
is a Kleinian group, consider the graph with one vertex $w_i$ for each
isomorphism class $\{\W_i\mid 1\le i\le r\}$ of irreducible representations
of $\Gamma$, and $n_{ij}$ edges between vertices $w_i$ and $w_j$, where the
positive integers $n_{ij}$ are the Clebsch-Gordon coefficients
\begin{equation} \label{McKay}
n_{ij} = \dim_\C \Hom_\Gamma(\V \o \W_i , \W_j) .
\end{equation}
The resulting graph is the Dynkin diagram of an irreducible simply-laced
affine Lie algebra; equivalently, the graph is connected, the numbers
$n_{ij}$ are equal to $0$ or $1$, and the Cartan matrix defined by
$$
A_{ij} = 2\delta_{ij} - n_{ij}
$$
is positive semi-definite, with one-dimensional null-space. In fact, the
null-space is spanned by the vector whose $i$th component is the dimension
of $\W_i$, since by \eqref{McKay},
$$
\sum_{j=1}^r (2\delta_{ij} - n_{ij}) \dim(\W_j) = 0 .
$$

\subsection*{Examples of the McKay correspondence}
\subsubsection*{Cyclic groups}
If $\Gamma=\CC_n$ is a cyclic group, let $\chi$ be the primitive character
characterized by $\chi(T)=\eps_n$. The irreducible representations of the
cyclic group $\CC_n$ are the powers $\{\chi^i\mid 0\le i<n\}$ of
$\chi$. Since $\V\o\chi^i\cong\chi^{i+1}\oplus\chi^{i-1}$, the associated
graph is a circuit with $n$ vertices: the Dynkin diagram $\hat{A}_{n-1}$.

\subsubsection*{Quaternionic groups}
The McKay correspondence associates to the quaternionic group $\QQ_{4n}$
the Dynkin diagram of $\hat{D}_{n-1}$:
$$
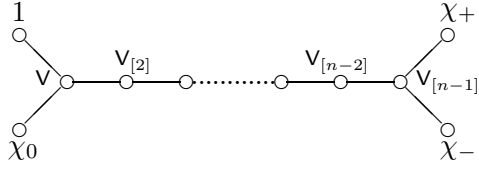
\begin{picture}(556,240)(141,-220)
\thinlines
\put(150,-36){\circle{15}}
\put(150,-156){\circle{15}}
\put(210,-96){\circle{15}}
\put(285,-96){\circle{15}}
\put(360,-96){\circle{15}}
\put(480,-96){\circle{15}}
\put(555,-96){\circle{15}}
\put(630,-96){\circle{15}}
\put(690,-36){\circle{15}}
\put(690,-156){\circle{15}}
\put(157,-43){\line( 1,-1){45}}
\put(157,-148){\line( 1, 1){45}}
\put(217,-96){\line( 1, 0){60}}
\put(292,-96){\line( 1, 0){60}}
\put(487,-96){\line( 1, 0){60}}
\put(562,-96){\line( 1, 0){60}}
\put(637,-88){\line( 1, 1){45}}
\put(637,-103){\line( 1,-1){45}}
\multiput(370,-102)(8.97436,0.00000){12}{\makebox(1.6667,11.6667){.}}
\put(140,-19){$1$}
\put(135,-186){$\chi_0$}
\put(680,-12){$\chi_+$}
\put(680,-186){$\chi_-$}
\put(170,-100){$\scriptstyle\V$}
\put(270,-75){$\scriptstyle\V_{[2]}$}
\put(510,-75){$\scriptstyle\V_{[n-2]}$}
\put(650,-100){$\scriptstyle\V_{[n-1]}$}
\end{picture}$$
The irreducible representations $\V_{[i]}$, $1\le i\le n-1$, are
two-dimensional, and $\V\cong\V_{[1]}$. The group $\QQ_{4n}$ has four
one-dimensional characters, as follows:
$$\begin{tabular}{|C|R|R|R|} \hline
\rho & \rho(S) & \rho(T) & \rho(U) \\ \hline
1 & 1 & 1 & 1 \\
\chi_0 & 1 & -1 & -1 \\
\chi_+ & -1 & i^n & - i^n \\
\chi_- & -1 & - i^n & i^n \\ \hline
\end{tabular}$$
Note that $\chi_+\V_{[i]}\cong\chi_-\V_{[i]}\cong\V_{[n-i]}$ and that
$\chi_0\V_{[i]}\cong\V_{[i]}$.

\subsubsection*{The binary octahedral group}
The cases $\TT$, $\OO$ and $\II$ of the McKay correspondence correspond
respectively to the affine Dynkin diagrams $\hat{E}_6$, $\hat{E}_7$ and
$\hat{E}_8$. Of these, we only need the case of $\OO$ in this paper; its
Dynkin diagram is as follows:
$$
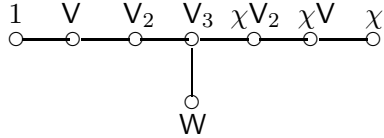
\begin{picture}(466,200)(51,-230)
\put(210,-91){\circle{15}}
\put(281,-91){\circle{15}}
\put(360,-91){\circle{15}}
\put(131,-91){\circle{15}}
\put(60,-91){\circle{15}}
\put(281,-170){\circle{15}}
\put(431,-91){\circle{15}}
\put(510,-91){\circle{15}}
\put(217,-91){\line( 1, 0){60}}
\put(292,-91){\line( 1, 0){60}}
\put(67,-91){\line( 1, 0){60}}
\put(142,-91){\line( 1, 0){60}}
\put(281,-103){\line( 0,-1){60}}
\put(367,-91){\line( 1, 0){60}}
\put(442,-91){\line( 1, 0){60}}
\put(50,-70){$1$}
\put(117,-70){$\V$}
\put(195,-70){$\V_2$}
\put(270,-70){$\V_3$}
\put(330,-70){$\chi\V_2$}
\put(415,-70){$\chi\V$}
\put(500,-70){$\chi$}
\put(265,-210){$\W$}
\end{picture}$$
The unique non-trivial one-dimensional character $\chi$ is characterized by
$\chi(S)=\chi(U)=-1$. Note that $\chi\V_3\cong\V_3$ and that
$\chi\W\cong\W$.

\section{The automorphism group of a hyperelliptic curve}

Denote by $\I_{2g+2}$ the affine variety of polynomials of degree $2g+2$
with non-vanishing discriminant. We identify a polynomial $f\in\I_{2g+2}$
with the binary form $y^{2g+2}f(x/y)$.

If $f\in\I_{2g+2}$, consider the affine varieties
$V\bigl(z^2-f(x)\bigr)\subset\Spec\C[x,z]$ and
$V\bigl(\tilde{z}^2-\tilde{x}^{2g+2}f(1/\tilde{x})\bigr)
\subset\Spec\C[\tilde{x},\tilde{z}]$. The hyperelliptic curve $C_f$
associated to $f$ is the smooth curve defined by gluing
$V\bigl(z^2-f(x)\bigr)$ and
$V\bigl(\tilde{z}^2-\tilde{x}^{2g+2}f(1/\tilde{x})\bigr)$ by the
identification $\bigl(\tilde{x},\tilde{z}\bigr)=\bigl(1/x,z/x^{g+1}\bigr)$.

The involution $\sigma:C_f\to C_f$ defined by $\sigma(x,z)=(x,-z)$ is
called the hyperelliptic involution of $C_f$; it acts on $H^0(C_f,\Om)$ by
$-\Id$ and is in the centre of the automorphism group $\Aut(C_f)$ of $C_f$.

\begin{lemma}
The curve $C_f$ has genus $g$, and $H^0(C_f,\Om)$ has basis
$\om_i=x^i\,dx/z$, $0\le i<g$.
\end{lemma}
\begin{proof}
The fixed points of the hyperelliptic involution are the $2g+2$ Weirstrass
points of $C_f$. The projection $(x,z)\mapsto x$ exhibits $C_f$ as a double
cover of $\PP^1$, ramified at the roots of $f$; thus, its genus equals $g$.

On the affine variety $V\bigl(z^2-f(x)\bigr)$, we have
$2z\,dz=f'(x)dx$; thus $\om_i = 2\,x^i\,dz/f'(x)$. Since the functions $z$
and $f'(x)$ have no common zeroes (the polynomial $f(x)$ has no multiple
roots), we conclude that the differentials $\om_i$ are regular on
$V\bigl(z^2-f(x)\bigr)$, so long as $i\ge0$.

Let $\tilde{f}(x)=x^{2g+2}f(1/x)$. On the affine variety
$V\bigl(\tilde{z}^2-\tilde{x}^{2g+2}f(1/\tilde{x})\bigr)
=V\bigl(\tilde{z}^2-\tilde{f}(\tilde{x})\bigr)$, we have
$2\tilde{z}\,d\tilde{z}=\tilde{f}'(\tilde{x})d\tilde{x}$. Since
$$
\om_i = - \tilde{x}^{g-i-1}\,d\tilde{x}/\tilde{z} = - 2 \,
\tilde{x}^{g-i-1} \, d\tilde{x} / \tilde{f}(\tilde{x}) ,
$$
the differentials $\om_i$ are regular on
$V\bigl(\tilde{z}^2-\tilde{x}^{2g+2}f(1/\tilde{x})\bigr)$ so long as $i<g$.

We have exhibited $g$ linearly independent algebraic one-forms on $C_f$;
since $C_f$ has genus $g$, they form a basis of $H^0(C_f,\Om)$.
\end{proof}

The group $\SL(2,\C)\times\C^\times$ acts by rational transformations on
$\Spec\C[x,z]$ by the formula
$$
\bigl( \bigl( \begin{smallmatrix} a&b \\ c&d \end{smallmatrix} \bigr) , u
\bigr) \* (x,z) = \bigl( \tfrac{ax+b}{cx+d} , \tfrac{uz}{(cx+d)^{g+1}}
\bigr) .
$$
If $f\in\I_{2g+2}$, the subgroup of elements of $\SL(2,\C)\times\C^\times$
which preserve the subvariety $V\bigl(z^2-f(x)\bigr)$ is a group of the
form
$$
\Gamma(\rho) = \bigl\{ ( \gamma,u) \mid \gamma \in \Gamma , u^2 =
\rho(\gamma) \bigr\} \subset \SL(2,\C) \times \C^\times ,
$$
where $\Gamma=\Gamma_f$ is the finite subgroup of $\SL(2,\C)$ consisting of
elements whose action on $\PP^1$ preserves the set of roots of $f$, and
$\rho=\rho_f$ is an even character of $\Gamma_f$. We have the short exact
sequence
$$
0 \to \bigl\<(-\Id,(-1)^{g+1})\bigr\> \to \Gamma(\rho) \to \Aut(C_f) \to 0
.
$$

Given a $\Gamma$-module $\W$ and an integer $n$, let $\W(n)$ be the
$\Gamma(\rho)$-module with underlying vector space $\W$ on which the
element $(\gamma,u)\in\Gamma(\rho)$ acts by
$$
(\gamma,u) \* w = u^n (\gamma\*w) .
$$
We have the isomorphisms $\W(n+2)\cong\rho\o\W(n)$ and
$\W(n)^\Dual\cong\W(-n)$. With this notation, the irreducible
representations of $\Aut(C_f)$ have the form $\W(n)$, where $\W$ is an
isomorphism class of irreducible representations of $\Gamma_f$, and
$n\equiv g\pmod{2}$ (respectively $n\equiv g+1\pmod{2}$) if $\W$ is even
(respectively odd).
\begin{proposition} \label{Schur}
As a representation of $\Aut(C_f)$, $H^0(C_f,\Om)\cong\V_{g-1}(-1)$.
\end{proposition}
\begin{proof}
Under the action of $(\gamma,u)\in\Gamma(\rho)$, $\om_i$ transforms into
$$
\bigl( \tfrac{ax+b}{cx+d} \bigr)^i \bigl( \tfrac{uz}{(cx+d)^{g+1}}
\bigr)^{-1} \tfrac{dx}{(cx+d)^2} = u^{-1} (ax+b)^i (cx+d)^{g-i-1} dx/z .
$$
Expanding the right-hand side in terms of the basis $\om_i$, we recover the
action of $\Aut(C_f)$ on $\V_{g-1}(-1)$.
\end{proof}

\begin{corollary} \label{basic}
As a representation of $\Aut(C_f)$,
$H^1(C_f,\C)\cong\V_{g-1}(1)\oplus\V_{g-1}(-1)$.
\end{corollary}
\begin{proof}
$H^1(C_f,\C)\cong H^0(C_f,\Om)\oplus H^1(C_f,\CO) \cong H^0(C_f,\Om)\oplus
H^0(C_f,\Om)^\Dual$
\end{proof}

Given a Kleinian group $\Gamma$ and a character $\rho$, let
$\I_{2g+2}(\Gamma,\rho)$ be the subvariety of $\I_{2g+2}$ consisting of
polynomials such that the pair $(\Gamma_f,\rho_f)$ is conjugate to
$(\Gamma,\rho)$.

The quotient $\CH_g$ of $\I_{2g+2}$ by the group
$\bigl(\SL(2,\C)\times\C^\times\bigr)/\bigl\<(-\Id,(-1)^{g+1})\bigr\>$ is
the moduli space of hyperelliptic curves of genus $g$. It is a complex
orbifold of dimension $2g-1$, stratified by the images $\CH_g(\Gamma,\rho)$
of the subvarieties $\I_{2g+2}(\Gamma,\rho)$. It carries a local system
$\VV$ whose fibre at $[f]$ is isomorphic to $H^1(C_f,\C)$; by Corollary
\ref{basic}, this local system has underlying vector bundle
$$
\I_{2g+2} \times_{(\SL(2,\C)\times\C^\times)/\<(-\Id,(-1)^{g+1})\>} \bigl[
\V_{g-1}(1) \oplus \V_{g-1}(-1) \bigr] .
$$

The same argument which was used in the introduction to calculate
$e(\CM_{1,1},\Schur^k\VV)$ proves the following result. This proposition
will be used in Section~4 to calculate the Euler characteristics
$e_2(1^k2^\ell)$.
\begin{proposition} \label{Basic}
$$
e\bigl( \CH_g,\WW(1^{k_1}\dots g^{k_g}) \bigr) = \sum_{(\Gamma,\rho)}
e\bigl( \CH_g(\Gamma,\rho) \bigr) \* \textstyle \dim \Bigl( \bigotimes_{i=1}^g
\Schur^{k_i}\bigl(\Wedge^i\bigl(\V_{g-1}(1)\oplus\V_{g-1}(-1)\bigr)\bigr)
\Bigr)^{\Gamma(\rho)}
$$
\end{proposition}

\section{The stratification of $\CH_2$}

We now specialize to genus $2$; this case is special, in that $\CH_2$ is
identical with the moduli space $\CM_2$ of smooth projective curves of
genus $2$.

Bolza \cite{Bolza} has shown that the stratification
$\CH_2=\coprod\CH_2(\Gamma,\rho)$ has seven strata. In Figure~1, we give a
diagram showing these strata, as well as two more pieces of data which we
will need: the Euler characteristics $e(\CH_2(\Gamma,\rho))$ of the strata,
and a normal form for polynomials in $\I_{2g+2}(\Gamma,\rho)$. Since no two
distinct strata have the same isotropy group $\Gamma$, we may, without
ambiguity, denote the stratum $\CH_2(\Gamma,\rho)$ by $\CH_2(\Gamma)$.

\begin{figure}
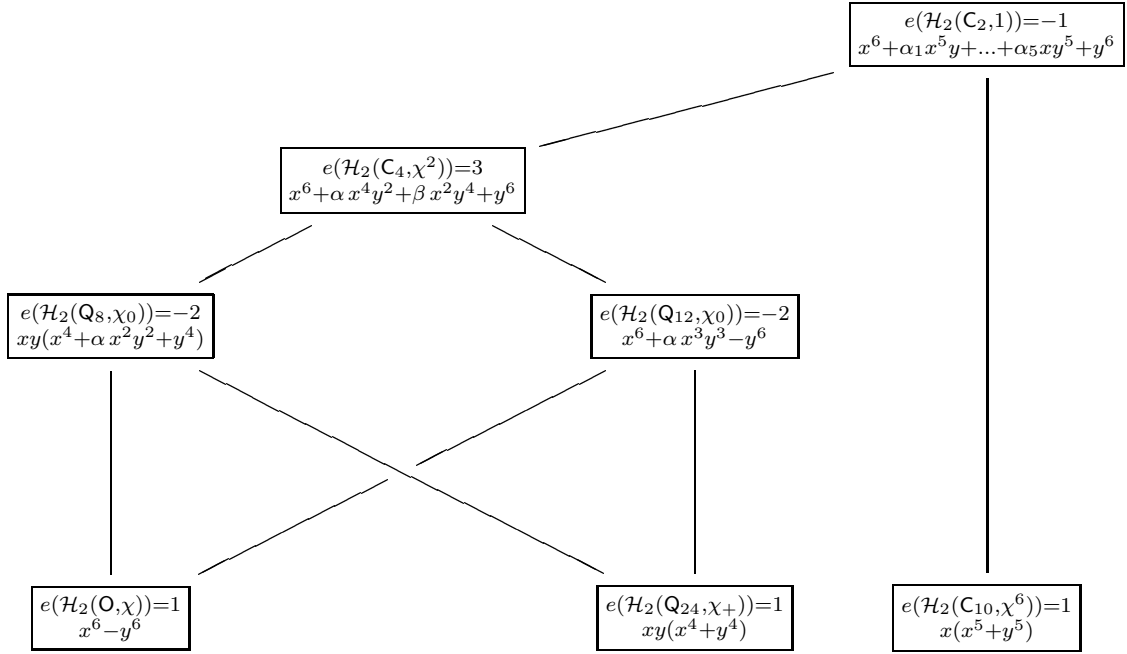

\caption{The Bolza stratification of $\CH_2$}
$$\begin{diagram}
\node[7]{\framebox{$\substack{e(\CH_2(\CC_2,1))=-1\\
x^6+\alpha_1x^5y+\ldots+\alpha_5xy^5+y^6}$}} \arrow[2]{wsw,-}
\arrow[8]{s,-} \\[2]
\node[3]{\framebox{$\substack{e(\CH_2(\CC_4,\chi^2))=3 \\
x^6+\alpha\,x^4y^2+\beta\,x^2y^4+y^6}$}} \arrow[2]{sw,-} \arrow[2]{se,-}
\\[2] \node{\framebox{$\substack{e(\CH_2(\QQ_8,\chi_0))=-2\\
xy(x^4+\alpha\,x^2y^2+y^4)}$}} \arrow[4]{s,-} \arrow[4]{se,-}
\node[4]{\framebox{$\substack{e(\CH_2(\QQ_{12},\chi_0))=-2\\
x^6+\alpha\,x^3y^3-y^6}$}} \arrow[4]{s,-} \arrow[2]{sw,-} \\[2] \node[3]{}
\arrow[2]{sw,-} \\[2]
\node{\framebox{$\substack{e(\CH_2(\OO,\chi))=1\\x^6-y^6}$}}
\node[4]{\framebox{$\substack{e(\CH_2(\QQ_{24},\chi_+))=1\\xy(x^4+y^4)}$}}
\node[2]{\framebox{$\substack{e(\CH_2(\CC_{10},\chi^6))=1\\x(x^5+y^5)}$}}
\end{diagram}$$
\end{figure}

In this section, we calculate the Euler characteristics of these
strata. Since $\CH_2(\CC_{10})$, $\CH_2(\QQ_{24})$ and $\CH_2(\OO)$ each
consist of precisely one point, it is clear that they have Euler
characteristic $1$. Since $\CH_2$ is contractible, it also has Euler
characteristic $1$; thus the Euler characteristics of all of the strata add
up to $1$. It remains to calculate the Euler characteristics
$e(\CH_2(\CC_4))$, $e(\CH_2(\QQ_8))$ and $e(\CH_2(\QQ_{12}))$. Of these,
the first is the most difficult; in calculating it, we use Clebsch's
classification of the covariants of binary sextics.

\subsubsection*{The Euler characteristic of $\CH_2(\CC_4)$}
If $f$ and $g$ are binary forms of degree $k$ and $\ell$ respectively,
define their $p$th Ueberschiebung $(f,g)_p$ by the formula
$$\textstyle (f,g)_p = \frac{(k+\ell-p)!}{(k+\ell)!} \bigl(
\frac{\p^2}{\p_x\p_\eta} - \frac{\p^2}{\p_y\p_\xi} \bigr)^p f(x,y)
g(\xi,\eta) \big|_{x=\xi,y=\eta} .
$$
The binary form $(f,g)_p$ is a joint covariant of $f$ and $g$.

For a proof of the following result, see Section 130 of Gordan
\cite{Gordan}.
\begin{lemma} \label{C(l,m)}
If $l$ and $m$ are a pair of quadratic forms, let $C(l,m)$ be the joint
invariant
$$
C(l,m) = \Bigl|\begin{smallmatrix} (l,l)_2 & (l,m)_2 \\ (m,l)_2 & (m,m)_2
\end{smallmatrix}\Bigr| .
$$
The quadratic forms $l$ and $m$ may be simultaneously diagonalized (i.e.\
there are coordinates $\xi$ and $\eta$ such that $l=a\,\xi^2+b\,\eta^2$ and
$m=c\,\xi^2+d\,\eta^2$) if and only if $C(l,m)\ne0$.
\end{lemma}

Given a binary sextic $f$, define a quartic covariant $i=(f,f)_4$ of degree
$2$, and quadratic covariants
$$
l=(i,f)_4 , \quad m=(i,l)_2 , \quad n=(i,n)_2
$$
of degrees $3$, $5$ and $7$ respectively. Let $R=-2\,((l,m)_1,n)_2$. Using
his symbol calculus for covariants of binary forms, Clebsch has shown
(\S113, \cite{Clebsch}) that
\begin{enumerate}
\item if $R\ne0$, $f$ is a cubic polynomial in the quadratic forms $l$, $m$
and $n$;
\item if $R=0$ and $C(l,m)\ne0$, $f$ is a cubic polynomial in the quadratic
forms $l$ and $m$.
\end{enumerate}
In each case, the coefficients of the representation are explicit rational
invariants of $f$. For an exposition of the proofs, see Gordan (\S29,
\cite{Gordan}).

We are interested in the second case above. By the condition $C(l,m)\ne0$
and Proposition~\ref{C(l,m)}, we see that there are coordinates $\xi$ and
$\eta$ such that $f\in\C[\xi^2,\eta^2]$. Furthermore, provided the
discriminant of $f$ is nonzero, we can rescale the coordinates $\xi$ and
$\eta$ in such a way that the coefficients of $\xi^6$ and $\eta^6$ equal
$1$.

In conclusion, a binary sextic $f$ with non-vanishing discriminant such
that $R=0$ and $C(l,m)\ne0$ is equivalent to a sextic in the normal form
\begin{equation} \label{alphabeta}
f(x,y) = x^6 + \alpha\,x^4y^2 + \beta\,x^2y^4 + y^6 .
\end{equation}
By Section I of Bolza \cite{Bolza}, these are precisely the sextics whose
image in $\CH_2$ lies in $\CH_2(\CC_4)$.

The discriminant of the normal form $f$ of \eqref{alphabeta} equals
$-64\,X^2$, where
$$
X = 4(\alpha^3+\beta^3) - \alpha^2\beta^2 - 18\alpha\beta + 27 ,
$$
while the invariant $(f,f)_6$ equals $\frac{2}{15}Y+2$, where
$$
Y = \alpha\beta .
$$
Conversely, from the functions $X$ and $Y$, we may recover the normal form
$f$ of \eqref{alphabeta} in the following way. Since
$$
16(\alpha^3-\beta^3)^2 = (X+Y^2+18\,Y-27)^2 - 64\,Y^3 ,
$$
we see that we can recover the divisor $(\alpha^3)+(\beta^3)$ in $\C$ from
$X$ and $Y$. From this, it is easy to recover the coefficients
$(\alpha,\beta)$, up to the action of the dihedral group generated by the
transformations $(\alpha,\beta)\mapsto(e^{2\pi i/3}\alpha,e^{-2\pi
i3}\beta)$ and $(\alpha,\beta)\mapsto(\beta,\alpha)$. Since these
transformations leave the functions $X$ and $Y$ invariant, we conclude that
$X$ and $Y$ are global coordinates on the stratum $\CH_2(\CC_4)$.

For the normal form $f$ of \eqref{alphabeta}, we have
$$
C(l,m) = - \frac{24^2}{15^{12}} \bigl( X+128(Y-9) \bigr)^2 \bigl(
(X+Y^2+18\,Y-27)^2-64\,Y^3 \bigr).
$$
Denoting the divisors $(X)$, $(X+128(Y-9))$ and
$((X+Y^2+18\,Y-27)^2-64\,Y^3)$ in $\C^2$ by $\Delta_0$, $\Delta_1$ and
$\Delta_2$, we conclude that
$$
\CH_2(\CC_4) \cong \C^2 \setminus \Delta_0 \cup \Delta_1 \cup \Delta_2 .
$$
In fact, Bolza shows that $\Delta_1$ corresponds to the stratum
$\CH_2(\QQ_8)$, while $\Delta_2$ corresponds to the stratum
$\CH_2(\QQ_{12})$.

It is now quite easy to calculate $e(\CH_2(\CC_4))$. Since $\Delta_0$ and
$\Delta_1$ are lines in $\C^2$, they have Euler characteristic $1$. The
projection $(X,Y)\mapsto Y$ displays $\Delta_2$ as a double cover of the
line ramified at $0$, showing that it too has Euler characteristic $1$. As
for the intersections of the three divisors, we have
$$
\Delta_0\cap\Delta_1 = \{(0,9)\} , \quad
\Delta_0\cap\Delta_2 = \{(0,9),(0,1)\} , \quad
\Delta_1\cap\Delta_2 = \{(0,9),(-2048,25),(-27648,225)\} .
$$
Combining these results, we conclude that $e(\CH_2(\CC_4))=3$.

\subsubsection*{The Euler characteristic of $\CH_2(\QQ_8)$}
Two binary sextics
$$
f(x,y) = xy(x^4+\alpha_ix^2y^2+y^4) , \quad i=1,2 ,
$$
are equivalent if and only if $\alpha_1^2=\alpha_2^2$; thus, we may
parametrize the stratum $\CH_2(\QQ_8)$ by $\alpha^2$. Three values of
$\alpha^2$ are excluded:
\begin{enumerate}
\item when $\alpha^2=4$, the sextic has vanishing discriminant;
\item when $\alpha^2=0$, the sextic has automorphism group $\OO$;
\item when $\alpha^2=100/9$, the sextic has automorphism group $\QQ_{24}$.
\end{enumerate}
It follows that the parameter $\alpha^2$ identifies the stratum
$\CH_2(\QQ_8)$ with $\C\setminus\bigl\{0,4,\tfrac{100}{9}\bigr\}$, and
hence that $e(\CH_2(\QQ_8))=-2$.

\subsubsection*{The Euler characteristic of $\CH_2(\QQ_{12})$}
Two binary sextics
$$
f(x,y) = x^6+\alpha_ix^3y^3-y^6 , \quad i=1,2 ,
$$
are equivalent if and only if $\alpha_1^2=\alpha_2^2$; thus, we may
parametrize the stratum $\CH_2(\QQ_{12})$ by $\alpha^2$. Three values are
excluded:
\begin{enumerate}
\item when $\alpha^2=-4$, the sextic has vanishing discriminant;
\item when $\alpha^2=0$, the sextic has automorphism group $\QQ_{24}$;
\item when $\alpha^2=50$, the sextic has automorphism group $\OO$.
\end{enumerate}
It follows that the parameter $\alpha^2$ identifies the stratum
$\CH_2(\QQ_{12})$ with $\C\setminus\{0,-4,50\}$, and hence that
$e(\CH_2(\QQ_{12}))=-2$.

\section{The calculation of $f_2(u,v)$}

We have seen in Proposition \ref{Basic} that
$$
f_2(u,v) = \sum_{(\Gamma,\rho)} e(\CH_2(\Gamma,\rho)) \*
\sum_{k,\ell=0}^\infty u^k v^\ell \dim \Bigl( \Schur^k\bigl(\VV\bigr) \o
\Schur^\ell\bigl(\Wedge^2\VV\bigr) \Bigr)^{\Gamma(\rho)} .
$$
where by Corollary \ref{basic}, $\VV=\V(1)\oplus\V(-1)$ and
$\Wedge^2\VV=\V_2\oplus\rho\oplus1\oplus\rho^{-1}$.
\begin{proposition}
\begin{multline*}
\sum_{k,\ell=0}^\infty u^kv^\ell \dim \bigl( \Schur^k\VV \o
\Schur^\ell\Wedge^2\VV \bigr)^{\Gamma(\rho)} 
\\ = \dim \Biggl( \frac{1+u^2(\V^2+\rho+\rho^{-1})+u^4}
{\bigl(1-u^2(\V^2-2)\rho+u^4\rho^2\bigr)
\bigl(1-u^2(\V^2-2)\rho^{-1}+u^4\rho^{-2}\bigr)} \\
\frac{1}{\bigl( 1 - v(\V^2-2) + v^2 \bigr) (1-v\rho) (1-v)^2
(1-v\rho^{-1})} \Biggr)^{\Gamma(\rho)} .
\end{multline*}
\end{proposition}
\begin{proof}
All of the following calculations are performed in the ring
$R(\Gamma)(u,v)$. Using the formulas $\sigma_t(A\oplus B) =
\sigma_t(A)\*\sigma_t(B)$ and $\sigma_t(A)=\lambda_{-t}(A)^{-1}$, we see
that
\begin{multline*}
\sum_{k,\ell=0}^\infty u^kv^\ell \Schur^k\VV \o \Schur^\ell\Wedge^2\VV 
= \sigma_u(\VV) \* \sigma_v(\Wedge^2\VV) \\
\begin{aligned}
{} &= \sigma_u(\V(1)) \* \sigma_u(\V(-1)) \* \sigma_v(\V_2) \*
\sigma_v(\rho) \* \sigma_v(1) \* \sigma_v(\rho^{-1}) \\
{} &= \lambda_{-u}(\V(1))^{-1} \* \lambda_{-u}(\V(-1))^{-1} \*
\lambda_{-v}(\V_2)^{-1} \* \lambda_{-v}(\rho)^{-1} \* \lambda_{-v}(1)^{-1}
\* \lambda_{-v}(\rho^{-1})^{-1} \\
{} &= \bigl(1-u\V(1)+u^2\rho\bigr)^{-1} \*
\bigl(1-u\V(-1)+u^2\rho^{-1}\bigr)^{-1} \* \\
{} & \qquad \bigl( 1 - v(\V^2-1) + v^2(\V^2-1) - v^3 \bigr)^{-1} \*
(1-v\rho)^{-1} \* (1-v)^{-1} \* (1-v\rho^{-1})^{-1} .
\end{aligned}
\end{multline*}
The third factor of the denominator may be simplified by the factorization
$$
\bigl( 1 - v(\V^2-1) + v^2(\V^2-1) - v^3 \bigr)^{-1}
= \bigl( 1 - v(\V^2-2) + v^2 \bigr)^{-1} (1 - v)^{-1} .
$$
To simplify the factors involving the variable $u$, we use the formulas
$$
1-u\V(1)+u^2\rho =
\frac{\bigl(1-u\V(1)+u^2\rho\bigr)\bigl(1+u\V(1)+u^2\rho\bigr)}
{1+u\V(1)+u^2\rho} = 
\frac{1-u^2(\V^2-2)\rho+u^4\rho^2} {1+u\V(1)+u^2\rho}
$$
and, similarly,
$$
1-u\V(-1)+u^2\rho^{-1} = \frac{1-u^2(\V^2-2)\rho^{-1}+u^4\rho^{-2}}
{1+u\V(-1)+u^2\rho^{-1}}.
$$
Hence
\begin{multline*}
\bigl(1-u\V(1)+u^2\rho\bigr)^{-1} \*
\bigl(1-u\V(-1)+u^2\rho^{-1}\bigr)^{-1} \\
\begin{aligned}
{} &= \frac{\bigl(1+u\V(1)+u^2\rho\bigr)\bigl(1+u\V(-1)+u^2\rho^{-1}\bigr)}
{\bigl(1-u^2(\V^2-2)\rho+u^4\rho^2\bigr)
\bigl(1-u^2(\V^2-2)\rho^{-1}+u^4\rho^{-2}\bigr)} \\
{} &= \frac{\bigl(1+u^2(\V^2+\rho+\rho^{-1}\bigr)+u^4\bigl)
+(u+u^3)\bigl(\V(1)\oplus\V(-1)\bigr)}
{\bigl(1-u^2(\V^2-2)\rho+u^4\rho^2\bigr)
\bigl(1-u^2(\V^2-2)\rho^{-1}+u^4\rho^{-2}\bigr)} .
\end{aligned}
\end{multline*}
No representation of $\Gamma(\rho)$ of the form $\W(n)$ with $n$ odd can
have a non-trivial space of invariants; we conclude that we may discard the 
terms which are odd in $u$ before taking the space of invariants under the
group $\Gamma(\rho)$.
\end{proof}

To apply this formula, we substitute for $\V$ the matrix $(n_{ij})$ with
entries $0$ and $1$ associated to the Dynkin diagram corresponding to
$\Gamma$ in the McKay correspondence, and for $\rho$ the permutation matrix
$(p_{ij})$ with entries
$$
p_{ij} = \dim_\C\Hom_\Gamma(\rho\o\W_i,\W_j) .
$$
In this way, we obtain a matrix with entries in $\Z(u,v)$; the desired
power series is the diagonal entry corresponding to the trivial
representation. We list the results for the seven cases of $(\Gamma,\rho)$
in Table~1, expressed as a sum of terms of the form $r(u,v)/s(u)t(v)$. We
have also listed the sum of the contributions for all strata other than
$\CH_2(\CC_{10})$; to obtain the formula for $f_2(u,v)$, we simply add this
total to the contribution for $\CC_{10}$.

\section{The equivariant Euler characteristic of $\CM_{2,n}$}

The following Euler characteristics are immediate from Table~1:
$$\begin{tabular}{|C||C|C|C|C|C|C|C|C|C|C|} \hline
(k,\ell) & (0,0) & (2,0) & (0,1) & (4,0) & (2,1) & (0,2) & (6,0) & (4,1) &
(2,2) & (0,3) \\ \hline
e_2(1^k2^\ell) & 1 & 0 & 0 & 0 & -1 & 0 & -1 & -1 & -1 & -3 \\ \hline
\end{tabular}$$
One might imagine from these data that all of the Euler characteristics
$e_2(1^k2^\ell)$ are negative: however, $e_2(1^{10})=1$.

Using the Leray-Serre spectral sequence, the $\SS_n$-equivariant Euler
characteristics of the moduli spaces $\CM_{g,n}$ may be expressed in terms
of the Euler characteristics $e_2(1^{k_1}\dots g^{k_g})$: by \cite{config},
we have
\begin{equation} \label{config}
\sum_{n=0}^\infty e_{\SS_n}(\CM_{g,n}) = e \biggl( (1+p_1)^2
\prod_{k=1}^\infty (1+p_k)^{-\frac{1}{k}\sum_{d|k}d!\,\mu(k/d)\ch_d(\VV)}
\biggr) .
\end{equation}
In this formula, we identify the virtual representation ring $R(\SS_n)$ the
symmetric group $\SS_n$ with the space of symmetric functions of degree
$n$, which, when tensored with $\Q$, is in turn isomorphic to the algebra
of polynomials of the power sums $p_k$.

In applying \eqref{config} in genus $2$, we may take advantage of the fact
that $e_2(1^k2^\ell)$ vanishes if $k$ is odd. We obtain the results listed
in the second column of Table~2. Substituting the values for
$e_2(1^k2^\ell)$, we obtain the equivariant Euler characteristics of
$\CM_{2,n}$, $0\le n\le7$. The dimensions of these virtual representations
of $\SS_n$ are listed in the third column: these numbers agree with the
Euler characteristics of $\CM_{2,n}$ calculated by Bini et al.\ \cite{BGP}.

\begin{table}
\caption{Contributions of the strata $\CH_2(\Gamma,\rho)$ to $f_2(u,v)$}
\begin{sideways}
$$
\begin{tabular}{|C|S|S|S|} \hline
(\Gamma,\rho) & \textstyle r(u,v) & \textstyle s(u) & \textstyle t(v)
\\[3pt] \hline\hline
(\CC_2,1) & u^4+6u^2+1 & (1-u^2)^4 & (1-v)^6 \\ \hline
(\CC_4,\chi^2) & (u^2+1)^2(v^4+6v^2+1)+16u^2(v^3+v) & (1-u^2)^4 &
(1-v)^2(1-v^2)^4 \\ \hline
(\QQ_8,\chi_0) & (u^8+u^6+4u^4+u^2+1)(v^4+1) +
(4u^8+14u^6+12u^4+14u^2+4)v^2 & (1-u^2)^2(1-u^4)^2 & (1-v)^2(1-v^2)^4 \\
& -(u^4-10u^2+1)(v^3+v) & (1-u^2)^4 & \\ \hline
(\QQ_{12},\chi_0) & (u^{12}+u^{10}+u^8+6u^6+u^4+u^2+1)(v^6+1) &
(1-u^2)^2(1-u^6)^2 & (1-v)(1-v^2)^4(1-v^3) \\
& +(3u^{12}+15u^{10}+31u^8+34u^6+31u^4+15u^2+3)(v^4+v^2) & & \\
& +2u^2 \bigl( (3u^8+5u^6+14u^4+5u^2+3)(v^5+v)
+2(u^4+1)(5u^4+11u^2+5)v^3 \bigr) & & \\ \hline
(\QQ_{24},\chi_+) &
(u^{16}+u^{14}+2u^{12}+4u^{10}+8u^8+4u^6+2u^4+u^2+1)(v^8+1) &
(1-u^4)^2(1-u^6)^2 & (1-v)^2(1-v^2)^3(1-v^6) \\
& +(2u^{16}+11u^{14}+24u^{12}+32u^{10}+30u^8+32u^6+24u^4+11u^2+2)(v^6+v^2)
& & \\
& -(u^{12}-3u^{10}-4u^8-12u^6-4u^4-3u^2+1)(v^7+v) & (1-u^2)^2(1-u^6)^2 & \\
& -(2u^{12}-5u^{10}-12u^8-18u^6-12u^4-5u^2+2)(v^5+v^3) & & \\
& +2(2u^{12}+2u^{10}+5u^8+6u^6+5u^4+2u^2+2)v^4 & & \\ \hline
(\OO,\chi) & (u^{20}+u^{18}+2u^{14}+6u^{12}+4u^{10}+6u^8+2u^6+u^2+1)(v^8+1)
& (1-u^4)(1-u^6)^2(1-u^8) & (1-v^2)^3(1-v^3)(1-v^4) \\
& +u^2(u^{16}+3u^{14}+15u^{12}+28u^{10}+26u^8+28u^6+15u^4+3u^2+1)(v^7+v) & & \\
& -
(u^{20}-5u^{18}-27u^{16}-57u^{14}-87u^{12}-106u^{10}-87u^8-57u^6-27u^4-5u^2+1)
(v^5+v^3) & & \\
& +2
(u^{20}+6u^{18}+18u^{16}+33u^{14}+46u^{12}+56u^{10}+46u^8+33u^6+18u^4+6u^2+1)
v^4 & & \\
& (u^8+3u^6+3u^2+1)(v^6+v^2) & (1-u^2)^4(1+u^4) & \\ \hline\hline
\text{Total for}
& u^2(2u^{16}+7u^{14}+15u^{12}+24u^{10}+24u^8+24u^6+15u^4+7u^2+2)(v^{10}+1)
& (1-u^4)(1-u^6)^2(1-u^8) & (1-v)^2(1-v^2)^2(1-v^4)(1-v^6) \\
\text{above strata}
& (2u^{16}+5u^{14}+15u^{12}+24u^{10}+28u^8+24u^6+15u^4+5u^2+2)(v^9+v)
& (1-u^2)^2(1-u^6)^2(1+u^4) & \\
& (u^{12}+13u^{10}+28u^8+36u^6+28u^4+13u^2+1)(v^8+v^2) & (1-u^2)^2(1-u^6)^2 &
\\
& +2(2u^{12}+14u^{10}+27u^8+34u^6+27u^4+14u^2+2)v^5 & & \\
& (3u^8+15u^6+4u^4+15u^2+3)(v^7+v^3) & (1-u^2)^4(1+u^4) & \\
& (2u^{12}+25u^{10}+47u^8+52u^6+47u^4+25u^2+2)(v^6+v^4) &
(1-u^2)^2(1-u^4)(1-u^8) & \\
\hline\hline
(\CC_{10},\chi^6) & (u^{12}-u^{10}+4u^8+4u^4-u^2+1)(v^4+1) &
(1-u^2)^3(1-u^{10}) & (1-v)^5(1-v^5) \\
& -(3u^{12}-11u^{10}+8u^8-8u^6+8u^4-11u^2+3)(v^3+v) & & \\
& +(5u^{12}-13u^{10}+16u^8-8u^6+16u^4-13u^2+5)v^2 & & \\ \hline
\end{tabular}
$$
\end{sideways}
\end{table}

\begin{table}
\caption{Calculation of $e_{\SS_n}(\CM_{2,n})$}
$$\begin{tabular}{|C|S|C|C|} \hline
n & & e_{\SS_n}(\CM_{2,n}) & e(\CM_{2,n}) \\ \hline\hline
0 & e_2 & 1 & 1 \\ \hline
1 & 2\,e_2\,s_1 & 2\,s_1 & 2 \\ \hline
2 & ( e_2 + e_2(2) ) s_2 + ( e_2 + e_2(1^2) ) s_{1^2}
& s_2 + s_{1^2} & 2 \\ \hline
3 & ( e_2(2) - e_2(1^2) ) s_3 + ( e_2(2) + e_2(1^2)
) s_{21} + 2\,e_2(1^2)\,s_{1^3} & 0 & 0 \\ \hline
4 & ( e_2 - e_2(1^2) - e_2(2) ) s_4 + ( - e_2 - e_2(1^2) +
e_2(2) ) s_{31} & s_4 - s_{31} & -4 \\
& + ( - e_2 + e_2(1^2) + e_2(2^2) ) s_{2^2} + ( e_2 - e_2(2) + e_2(1^22) )
s_{21^2} + ( e_2(1^2) + e_2(1^4) ) s_{1^4} & - s_{2^2} & \\ \hline
5 & ( 2\,e_2 - 2\,e_2(2) ) s_5 + ( e_2 - e_2(1^22) -
e_2(2^2) ) s_{41} & 2(s_5 + s_{41}) & 0 \\
& + ( - 3\,e_2 + 2\,e_2(1^2) + 2\,e_2(2) - e_2(1^22) + e_2(2^2) )
s_{32} & - 2\,s_{32} & \\
& + ( - e_2(1^2) + e_2(2) - e_2(1^4) - e_2(2^2) ) s_{31^2} & & \\
& + ( e_2 - 2\,e_2(2) + e_2(1^22) + e_2(2^2) ) s_{2^21} & & \\
& + ( e_2 - e_2(1^2) - e_2(2) + e_2(1^4) + e_2(1^22) )
s_{21^3} + 2\,e_2(1^4)\,s_{1^5} & & \\ \hline
6 & ( e_2 - e_2(1^2) - 2\,e_2(2) + e_2(1^22) ) s_6 + (
2\,e_2 - e_2(1^2) - e_2(2) + e_2(1^4) ) s_{51} & 2\,s_{51} & -24 \\
& + ( 2\,e_2(2) - e_2(2^2) ) s_{42} + ( 2\,e_2(1^2) + e_2(2)
- e_2(1^4) - 2\,e_2(1^22) - e_2(2^2) ) s_{41^2} & 2\,s_{41^2} & \\
& + ( -2\,e_2 + 2\,e_2(1^2) + e_2(1^4) + 2\,e_2(2^2) ) s_{3^2} &  -
2\,s_{3^2} & \\
& + ( - e_2 - e_2(1^2) - 2\,e_2(2) + e_2(1^22) + e_2(2^2) )
s_{321} & - 2\,s_{321} & \\
& + ( - 2\,e_2(1^2) + e_2(2) - e_2(1^4) - e_2(2^2) ) s_{31^3} & & \\
& + ( e_2 - e_2(1^2) - e_2(2) + e_2(1^4) + e_2(2^2) ) s_{2^3} & - 3\,
s_{2^3} & \\
& + ( e_2 - e_2(1^2) - e_2(2) + e_2(1^4) + e_2(1^22^2) )
s_{2^21^2} & s_{2^21^2} & \\
& + ( - e_2 + 2\,e_2(1^2) + e_2(2) - e_2(1^22) + e_2(1^42) )
s_{21^4} + ( e_2(1^6) + e_2(1^4) ) s_{1^6} & - ( s_{21^4} + s_{1^6} ) & \\
\hline
7 & (-e_2(1^4)+e_2(1^22)-e_2)s_7 +
(e_2(1^4)+2\,e_2(1^22)+e_2(2^2)-3\,e_2(1^2)-3\,e_2(2))s_{61} & - 2(s_7 +
s_{61} ) & 168 \\
& +(-2\,e_2(2^2)-2\,e_2(1^2)+2\,e_2(2)+2\,e_2)s_{52} & 2\,s_{52} & \\
& + (e_2(1^4)-2\,e_2(1^22)+e_2(2^2)+3\,e_2(1^2)-e_2(2))s_{51^2} &
2\,s_{51^2} & \\
& + (e_2(1^22)+e_2(2^2)+e_2(1^2)+e_2(2)-e_2)s_{43} & - 2\,s_{43} & \\
& +
(-e_2(1^22^2)-e_2(2^3)-2\,e_2(1^4)-e_2(1^22)-e_2(2^2)+4\,e_2(1^2)+3\,e_2(2)-1)
s_{421} & 4\,s_{421} & \\
& +
(-e_2(1^42)-e_2(1^22^2)-e_2(1^4)-e_2(1^22)+3\,e_2(1^2)+e_2(2)-1)s_{41^3} &
2\,s_{41^3} & \\
& + (2\,e_2(1^4)+e_2(1^22)+3\,e_2(2^2)+e_2(1^2)-3\,e_2(2)-1)s_{3^21} & -
2\,s_{3^21} & \\
& + (-e_2(1^22^2)+e_2(2^3)-e_2(1^4)+e_2(1^22)-e_2(1^2)-2\,e_2(2)-1)s_{32^2}
& - 4\,s_{32^2} & \\
& + (-e_2(1^42)-e_2(2^3)+2\,e_2(1^4)+2\,e_2(1^22)-4\,e_2(1^2)-e_2(2)+2\,e_2)
s_{321^2} & 4\,s_{321^2} & \\
& + (-e_2(1^6)-e_2(1^22^2)-e_2(1^4)+e_2(1^22)-e_2(1^2)-e_2)s_{31^4} & & \\
& + (e_2(1^22^2)+e_2(2^3)-e_2(1^22)-e_2(2^2)+e_2(2)+1)s_{2^31} & -
2\,s_{2^31} & \\
& + (e_2(1^42)+e_2(1^22^2)-2\,e_2(1^22)+3\,e_2(1^2)+e_2(2))s_{2^21^3} & & \\
& + (e_2(1^6)+e_2(1^42)-e_2(1^4)-e_2(1^22)+2\,e_2(1^2)+e_2(2)-e_2)s_{21^5}
+2\,e_2(1^6)s_{1^7} & - 2(s_{21^5} + s_{1^7}) & \\ \hline
\end{tabular}$$
\end{table}

\end{document}